\documentclass[
10pt, % Main document font size
a4paper, % Paper type, use 'letterpaper' for US Letter paper
oneside, % One page layout (no page indentation)
headinclude,footinclude, % Extra spacing for the header and footer
BCOR5mm, % Binding correction
]{scrartcl}

\usepackage[
nochapters, % Turn off chapters since this is an article        
beramono, % Use the Bera Mono font for monospaced text (\texttt)
eulermath,% Use the Euler font for mathematics
pdfspacing, % Makes use of pdftex’ letter spacing capabilities via the microtype package
dottedtoc % Dotted lines leading to the page numbers in the table of contents
]{classicthesis} % The layout is based on the Classic Thesis style

\usepackage{arsclassica} % Modifies the Classic Thesis package

\usepackage[T1]{fontenc} % Use 8-bit encoding that has 256 glyphs

\usepackage[utf8]{inputenc} % Required for including letters with accents

\usepackage{graphicx} % Required for including images
\graphicspath{{Figures/}} % Set the default folder for images

\usepackage{enumitem} % Required for manipulating the whitespace between and within lists

\usepackage{lipsum} % Used for inserting dummy 'Lorem ipsum' text into the template

\usepackage{subfig} % Required for creating figures with multiple parts (subfigures)

\usepackage{amsmath,amssymb,amsthm} % For including math equations, theorems, symbols, etc

\usepackage{varioref} % More descriptive referencing

\usepackage[numbers]{natbib} % So that \citep is recognized

\theoremstyle{definition} % Define theorem styles here based on the definition style (used for definitions and examples)

\theoremstyle{plain} % Define theorem styles here based on the plain style (used for theorems, lemmas, propositions)

\theoremstyle{remark} % Define theorem styles here based on the remark style (used for remarks and notes)

\hypersetup{
colorlinks=true, breaklinks=true, bookmarks=true,bookmarksnumbered,
urlcolor=webbrown, linkcolor=RoyalBlue, citecolor=webgreen, % Link colors
pdftitle={}, % PDF title
pdfauthor={\textcopyright}, % PDF Author
pdfsubject={}, % PDF Subject
pdfkeywords={}, % PDF Keywords
pdfcreator={pdfLaTeX}, % PDF Creator
pdfproducer={LaTeX with hyperref and ClassicThesis} % PDF producer
} % Include the structure.tex file which specified the document structure and layout

\usepackage{hyperref} % Links to footnotes
\usepackage{footnotebackref}

\makeatletter
\let\old@makefnmark\@makefnmark
\def\@makefnmark{\hbox{\@textsuperscript{\normalfont\hyper@linkstart{link}{footnote@\the\c@footnote}\old@makefnmark\hyper@linkend}}}
\makeatother

\title{\normalfont\spacedallcaps{Pascual Jordan's "Erweiterte Gravitationstheorie" - A Historical Analysis of its Mathematical Framework}} % The article title

\author{\spacedlowsmallcaps{Bernadette Lessel}} % The article author(s) - author affiliations need to be specified in the AUTHOR AFFILIATIONS block

\date{\today} % An optional date to appear under the author(s)

\begin{document}
	
%----------------------------------------------------------------------------------------
%	HEADERS
%----------------------------------------------------------------------------------------
	
	\renewcommand{\sectionmark}[1]{\markright{\spacedlowsmallcaps{#1}}} % The header for all pages (oneside) or for even pages (twoside)
	\lehead{\mbox{\llap{\small\thepage\kern1em\color{halfgray} \vline}\color{halfgray}\hspace{0.5em}\rightmark\hfil}} % The header style
	
	\pagestyle{scrheadings} % Enable the headers specified in this block
	
	%----------------------------------------------------------------------------------------
	%	TABLE OF CONTENTS & LISTS OF FIGURES AND TABLES
	%----------------------------------------------------------------------------------------
	
	\maketitle % Print the title/author/date block
	
	\setcounter{tocdepth}{2} % Set the depth of the table of contents to show sections and subsections only
	
	\tableofcontents % Print the table of contents

	%----------------------------------------------------------------------------------------
	%	ABSTRACT
	%----------------------------------------------------------------------------------------
	
\section*{Abstract}

This paper aims to highlight Pascual Jordan’s axiomatic definition of the covariant derivative, as set out in his 1952 textbook "Schwerkraft und Weltall". Developed in light of his \emph{Erweiterte Gravitationstheorie} — a projective reformulation of relativity theory that incorporates a variable gravitational constant — Jordan's definition resembles those in contemporary usage. The paper contextualises Jordan's work within the broader historical frameworks of differential geometry and projective relativity, with a particular focus on the Princeton relativity group led by Oswald Veblen and Luther Pfahl Eisenhart. It also provides a summary of Jordan's formalism, focusing particularly on his definition of the covariant derivative, as well as a brief history of the origin and development of the covariant derivative.

	%----------------------------------------------------------------------------------------
	%	AUTHOR AFFILIATIONS
	
%	\let\thefootnote\relax\footnotetext{* \textit{Universität Bonn}}
	
	%----------------------------------------------------------------------------------------
	
	\newpage % Start the article content on the second page, remove this if you have a longer abstract that goes onto the second page
	
%----------------------------------------------------------------------------------------
%	INTRODUCTION
%----------------------------------------------------------------------------------------
	
\section{Introduction}
%---------------------

The original idea for this article was for it to be a short note highlighting Pascual Jordan's innovative definition of the covariant derivative in his 1952 textbook "Schwerkraft und Weltall" \citep{Jordan1952}. In this book, Jordan gives an introduction to general relativity leading up to a treatment of projective relativity, on which basis Jordan developes his "extended theory of gravitation" (\emph{Erweiterte Gravitationstheorie}) including a non-constant gravitational "constant".\footnote{Jordan's \emph{Erweiterte Gravitationstheorie} is a mathematical treatment of an hypothesis originally due to Paul Dirac. In 1937, Dirac speculated that "the gravitational 'constant' must decrease with time, proportionally to $t^{-1}$" \citep{Dirac1937}. The conclusions, section \ref{s.reassesing}, will say a bit more about the physical background of Jordan's theory.} Jordan's covariant derivative therein is remarkable, because it is given in an operational, axiomatic way, which is completely uncommon at the time and which resembles very much the contemporary definition of a connection, as an operator on sections of vector bundles obeying a list of properties. And the reason Jordan comes to introduce the covariant derivative in this way, and this was the second point I wanted to highlight, is not an independent, intrinsically mathematical will to axiomatize, but his need to have a notion of a covariant derivative which is straighforwardly applicable to the projective case in a way that would naturally give way to his extension. In different words, it is the work on his personal, extended unified field theory that made him search for a definition which makes the development of the framework more straightforward and computations easier.

Researching other, earlier definitions of the covariant derivative brought about two obstructions to this original plan: A) Jordan was actually not the first person to come up with an axiomatic definition of the covariant derivative. There was (at least) Walther Mayer who in his textbook on differential geometry from 1930 gives a very similar definition \citep[p.156]{Mayer1930}. Mayer's motivation was analogous in that he wanted to abstract the necessary properties for application to a more general case. However, his was that of tensors defined on subspace of $\mathbb{R}^n$. It can reasonably be assumed that Jordan did not know about Mayer's textbook, or that, at least, he did not consciously use it. But still, he was not the first in a strict manner. And B) I came to realize that the history of projective theories of relativity (and unified field theory more generally) is very tightly connected to the advancement of geometry in the 1920's, including the understanding of the covariant derivative and the parallel transport therein. So, I thought, there is no way around incorporating this in my "note".

The aim of the paper has thus been extended to give a more general historical account of projective relativity. I begin in section \ref{s.formalism} with an introduction to the mathematical formalism of Jordan's extended theory of gravitation. This seems worthwile, because Jordan's publications of his theory were entirely in German, wherfore an English summary of the main lines of theory construction may enable a first contact with the topic for a non-German speaker. But also, it does not seem possible to showcase the remarkability of his covariant definition without such an account. 

In the intermediate section \ref{s.briefhistory}, I try to give a brief account of where the covariant definition historically originates from - and how it continued to develop until the early 1920, generating the concepts of parallel transport and connection along the way. In this section, I gratefully and mainly draw from already existing historiographical literature, especially from Karin Reich's excellent work on the topic.

In section \ref{s.projective}, I outline the development of one of the strands that lead to a projective reformulation of five-dimensional unified field theory, projective relativity, at the end of the 1920's. This strand is that of the Princeton relativity group around Oswald Veblen and Luther Pfahl Eisenhart. I am focusing on the work of this group, because even though Jordan's projective theory seems to owe more to the other strands, it is Veblen whom he cites mostly. Also, the approach of that group towards a generalized notion of geometry seems to demonstrate most exemplarily the struggles of geometry at that time. 
The final section, section \ref{s.reassesing}, then concludes by contextualizing Jordan's work on projective relativity in the light of its history. 

The formal and historical investigations in this text have mainly been carried out in order to bring to light the definition and usage of the covariant derivative in the period treated. For this, it is inevitable to having to work through some technical details of mathematical formalism. At the same time, I mostly decided not to go further in the formalism than the definition of the covariant derivative, being the fundamental process in these theories, and thus entirely neglected derived notions such as the curvature tensor. This is a deliberate choice, in order to put a limit somewhere. 

In any case, this paper should be understood as a genuine effort to integrate the histories of maths and physics, showing how the aim for a specific application within physics brought about development in math and vice versa. It is thus, hopefully, not only engaging for historians of maths interested in (differential) geometry, but also for historians of physics interested in general relativity and unified field theory. Jordan's textbook, in any case, laid the foundation for a new generation of relativists in Germany, most importantly his hown relativity group in Hamburg, whose first contact with research topics in that field was in projective relativity.
	
%----------------------------------------------------------------------------------------
%	METHODS
%----------------------------------------------------------------------------------------
	
\section{Jordan's "Erweiterte Gravitationstheorie" - The Formalism}\label{s.formalism}
%------------------------------------------------------------------
%------------------------------------------------------------------

This section summarizes the mathematical framework of Jordan's "extended theory of gravity" as presented in his textbook "Schwerkraft und Weltall" from 1952 \citep{Jordan1952}.\footnote{"Schwerkraft und Weltall" saw a second edition in 1955, which however did not differ at all with respect to what concern us here. (In other aspects it does differ, though.)} This book is both a consolidation of previous efforts on this theory, by Jordan himself and by his co-workers\footnote{In the 1952 edition of the book, Jordan mentions the names Fricke, Greßmann, Heckmann, Ludwig, Müller, Lüders. For the second edition Engelbert Schücking playes a major role.}, as well as a self-explanatory treatment of the material for late entrants. It includes a full treatment of Riemannian geometry and of Einstein's general relativity. In this respect the book could be used, and in fact had been used\footnote{In the preface to the second edition, Jordan states "Ich freue mich, aus mancherlei freundlichen Zuschriften ersehen zu haben, daß sich das Buch gerade in diesem Sinne bewährt hat."}, as an introductory text to general relativity, even if one was not interested in Jordan's theoretical extension.

It was imperative to Jordan, so he stated in the introduction of the book, to simplify the formal development of the "Riemann-Einsteinian theory", to be able to exectute the extension of the five-dimensional projective space  he had in mind \citep[p.IV]{Jordan1952}: 
	\begin{quote}
	However, from the outset, I present a development of the \emph{Riemann-Einstein} theory that has been modified in some respects and (I believe) simplified. The need to think intensively about possibilities for simplification arose in the course of investigations into five-dimensional projective relativity theory; the extension I have given to this theory would have been impossible without a prior radical simplification of the mathematical proofs of four-dimensional theory.
%	Ich gebe aber von Anfang an eine in mancher Hinsicht veränderte, und (wie ich glaube) vereinfachte Entwicklung der \emph{Riemann-Einstein}schen Theorie. Die Notwendigkeit, eindringlich über Möglichkeiten der Vereinfachung nachzudenken, ergab sich im Zuge der Untersuchungen zur fünfdimensional-projektiven Relativitätstheorie; die Erweiterung, welche ich dieser Theorie gegeben habe, wäre undurchführbar gewesen ohne eine vorherige radikale Vereinfachung der mathematischen Beweismittel der vierdimensionalen Theorie.
    \end{quote}
And one of the tools Jordan simplified, and which is the focus of this historical investigation, is the covariant derivative, of which Jordan was the first in the relativity literature to give an axiomatic definition, independent of the notion of a parallel transport.\footnote{It seems that in the entire literature, including mathematical literature, covering Riemannian geometry, only Walther Mayer was earlier in giving an axiomatic definition of the covariant derivative than Jordan. This was in his textbook series on differential geometry that he co-authored with Adalbert Duschek, \citep{Duschek1930}.} 

\subsection{Jordan's Axiomatic Definition of the Covariant Derivative}
%-----------------------------------------
%
Jordan approaches his axiomatic definition of the covariant derivative in a three-step procedure. First, he gives a constructive definition of the covariant derivative, in which the notion of a "coordinate system $x^{(k)}$ which is plane at point $P$"\footnote{"Ein im Punkte $P$ ebenes Koordinatensystem $x^{(k)}$"} plays a crucial role. Such a coordinate system is defined to have the following property: At point $P$ the metric components $g_{kl}$ in this coordinate system are constant in first order, i.e. $g_{kl|j}=0$.\footnote{Here and in the entire text we use original notation. The short-hand notation for the partial derivative of a tensor component, e.g.  $g_{kl|j}$, Jordan notes to have adapted from Pais.} With this, the definition for the covariant derivative reads the following \citep[p.31]{Jordan1952}: 

\emph{The tensor field $T^{...}_{...||k}$ is called covariant derivative of the tensor field $T^{...}_{...}$, if at every point $P$ and for every coordinate system which is plane at point $P$ the following property holds true:} 
\begin{equation}
	T^{l_1 l_2 ...l_s}_{h_1 h_2 ...h_r||k}=T^{l_1 l_2 ...l_s}_{k_1 k_2 ...k_r|k}.
\end{equation}
This definition is constructive, because the covariant derivative is constructed in a pointwise manner. For every point $P$ one needs to find a coordinate system which is plane at $P$ and then performs the ordinary partial derivative at that point and with respect to this coordinate system. The result is then transferred to arbitrary coordinate systems by means of known formula for coordinate transformations. 

That this definition is indeed well-defined of course depends on the existence of such plane coordinate systems at every point. That this is indeed so, is a two-page proof in Jordan's book. This proof is, at the same time, instructive as it introduces the Christoffel symbols $\Gamma^m_{jr}$ as components mediating the coordinate change via $x^{(k)}=\bar{x}^{(k)}-\frac{1}{2}\Gamma^k_{lj}\bar{x}^{(l)}\bar{x}^{(j)} + \text{higher orders}$, $\Gamma^k_{lj}=\Gamma^k_{jl}$. Here, $x^{(k)}$ denote the old and $\bar{x}^{(k)}$ denote the new coordinate functions and this equation is to be understood as a Taylor expansion where $P$ is at the origin of $x^{(k)}$. Turns out, if one chooses 
\begin{equation}
	\Gamma^m_{jr}=\frac{1}{2}g^{mh}\left(g_{hj|r}+g_{hr|j}-g_{jr|h}\right),
\end{equation}
then the new coordinate system $\bar{x}^{(k)}$ is plane in $P$. This equation is, of course, the known formula for Christoffel symbols.

As a second step, Jordan derives five consquences from the constructive definition for the covariant derivative:
\begin{enumerate}
	\item[I.] The usual rules for differentiating sums and products of tensors apply:
	      \begin{equation}\begin{aligned}
	      	\left(T^{...}_{...}+S^{...}_{...}\right)_{||k} &= T^{...}_{...||k}+S^{...}_{...||k}\\
	      	\left(T^{...}_{...}S^{...}_{...}\right)_{||k} &= T^{...}_{...||k}S^{...}_{...}+T^{...}_{...}S^{...}_{...||k}.
	      \end{aligned} \end{equation}
      
     \item[II.] The operation of covariant differentiation is interchangeable with the operation of contraction. That is, $T^{j...}_{j...||k}$ equals the contraction of $T^{j...}_{l...||k}$.
     
     \item[III.] The operation of covariant differentiation is interchangeable with the operation of raising or lowering indices: $g_{kl||j}=0$.
     
     \item[IV.] The gradient can be formed by covariant derivation: $\Phi_{||}=\Phi_{|}$.
     
     \item[V.] The rotation can be formed by covariant derivation: $a_{k||j}-a_{j||k}=a_{k|j}-a_{j|k}$.\footnote{Jordan identifies co- and contravariant vectors in a way that he does not distinguish between those after establishing the isomorphism of those spaces, thus also calls $a_k$ vectors.}
\end{enumerate}

Finally, Jordan gives an axiomatic definition of the covariant derivative, as the operation having the above five properties \citep[p.32]{Jordan1952}:
\begin{quote}
	The combination of these five facts I. to V. gives us an \emph{axiomatic characterization}\footnote{Every emphasis in this article is original.} of covariant differentiation: An operation that has these five properties is thus uniquely defined as covariant differentiation; all other properties of the covariant derivative can be derived from I. to V. without having to refer back to the definition in §8.
	%Die Zusammenstellung dieser fünf Tatsachen I. bis V. gibt uns eine \emph{axiomatische Kennzeichnung}\footnote{Every emphasis in this article is original.} der kovarianten Differentiation: Eine Operation, die diese fünf Eigenschaften hat, ist damit eindeutig als die kovariante Differentiation festgelegt; alle sonstigen Eigenschaften der kovarianten Ableitung können aus I. bis V. abgeleitet werden, ohne daß wir noch einmal auf die Definition von §8 zurückgreifen müssen.
\end{quote} 
The definition from §8 mentioned here is, of course, the constructive definition introduced above. Jordan justifies this axiomatic definition with its helpfulness in the five-dimensional projective theory: 
\begin{quote}
	This is a fact that will be of crucial help to us later in five-dimensional projective theory.
	%Das ist eine Tatsache, die uns später in der fünfdimensional-projektiven Theorie ganz entscheidende Hilfe geben wird. 
\end{quote}
%[Mention that Jordan reduces this definition to four axioms in Kundt lecture notes; perhaps one can reprint the pages in the article?]

Later, Jordan also gives a general formula to calculate the covariant derivative of a tensor of arbitrary type, despite asserting that "we hardly ever have to use it for calculations": 
\begin{equation}
		T^{l_1...l_s}_{k_1...k_r||h} = T^{l_1...l_s}_{k_1...k_r|h}
		                               +\sum_{i=1}^{s}\Gamma_{lh}^{l_i}T_{k_1...k_r}^{l_1...l_{i-1}ll_{i+1}...l_s}
		                               -\sum_{i=1}^r\Gamma_{k_ih}^k T^{l_1...l_s}_{k_1...k_{i-1}kk_{i+1}...k_r}.
\end{equation}
And it is only much after introducing the Riemannian curvature tensor $G$ (as the tensor satisfying $a_{k||l||j}-a_{k||j||l}=-G^m_{\cdot klj}a_m$)\citep[p.38]{Jordan1952} and deriving the typical, known formulas that Jordan discusses the notion of parallel transport - as part of a proof showing that from a vanishing Riemannian curvature tensor it follows that the metric is Euclidean \citep[pp.46--50]{Jordan1952}. In this proof, vectors are defined to be parallel transported along a curve as the unique solution of the differential equation of first order $a_{j||k}dx^k=0$ for the functions $a_j$ along a curve connecting two points, given an initial value (which can be considered the vector to be parallel transported along the curve). 
In any case, Jordan has thus demonstrated that the notion of parallel transport for him is really nothing but a helpful tool, certainly not fundamental.

\subsection{The Projective Structure of the Einstein-Maxwell System}
%---------------------------------------------------------------------------
As was mentioned in the introduction, Jordan emphasized that projective relativity is not about altering the Einstein-Maxwell theory, but about analyzing its mathematical structure. The first part of this analysis deals with the invariance group of the vacuum Einstein-Maxwell system \citep[§22]{Jordan1952}, the latter of which Jordan introduces to be the following set of Lagrange equations: 
\begin{equation}\begin{aligned}\label{eq.EM}
	G_{kl}+&\frac{\chi}{c^2}\left(F_{kj}F^{\cdot j}_l-\frac{1}{4}g_{kl}F_{hj}F^{hj}\right)=0;\\
		F^{kl}_{\ \ \ ||l} &=0,
	\end{aligned}
\end{equation}
resulting from the Lagrangian density $\left(G+\frac{\chi}{2c^2}F_{kl}F^{kl}\right)\sqrt{-g}$. Here, the $F_{kl}$ are defined as $F_{kl}=\Phi_{l|k}-\Phi_{k|l}$, $G_{kl}$ is the Ricci tensor and $\Phi_l$ is the electromagnetic four-potential. 

Jordan notes that the covariant equations \ref{eq.EM} are not only invariant against coordinate transformations, but also against gauge transformations, in which $\Phi_k$ is being substituted by 
\begin{equation}
	\Phi'_k=\Phi_k+\Phi_{|k},
\end{equation}
where $\Phi$ is an arbitrary (differentiable) function of the coordinates $x^{(k)}$. The group of gauge transformations is thus the additive group of arbitrary functions in $x^{(k)}$\footnote{Jordan writes, these functions are arbitrary, though obviously they need to be at least differentiable.}, which Jordan denotes with $\mathfrak{E}$. The group of fourdimensional coordinate transformations is denoted by $\mathfrak{K}$ and the system \ref{eq.EM} is invariant against the group that is jointly generated by $\mathfrak{E}$ and $\mathfrak{K}$. This generated total invariance group $\mathfrak{I}$, however, is not isomorphic to the direct product of $\mathfrak{E}$ and $\mathfrak{K}$. Instead, as Jordan shows, it is isomorphic to the group $\mathfrak{H}_5$ of all those transformations in five variables $X^0,X^1,...,X^4$, the new variables $\bar{X}^{\mu}$ of which are homogenous functions of degree one in the $X^{\nu}$. It is instrucitve to see that a general element of $\mathfrak{H}_5$ is of the form 
\begin{equation}
	\bar{X}^{\mu}=X^{\nu} F^{(\nu)}\left(X^1/X^0, X^2/X^0, X^3/X^0, X^4/X^0\right), 
\end{equation}
more important, however, is a differential property which is derived from the more general Euler theorem for homogenous functions\footnote{Jordan does not state the theorem (but obviously uses it), but it goes like this: For a continously differentiable funktion $f(x_1,...,x_n)$ which is homogenous of degree $N$, it holds $\sum_{i=1}^{n}\frac{\partial f}{\partial x_i} x_i=N f(x_1,...,x_n)$.}, 
\begin{equation}\label{eq.trafo}
	\bar{X}^{\mu}_{\ \ |\nu}X^\nu=\bar{X}^{\mu}.
\end{equation}

In this isomorphism, elements of $\mathfrak{E}$, denoted by $\left[\Phi(x^1,x^2,x^3,x^4)\right]$, and which constitute a normal subgroup within $\mathfrak{I}$, correspond to elements $\pm F\left(X^1/X^0, ..., X^4/X^0\right)=e^{\Phi\left(x^1,...x^4\right)}$, which in turn constitute a normal subgroup within $\mathfrak{H}_5$. Elements of $\mathfrak{K}$, which have the form $\bar{x}^k=f^k(x)$, correspond to elements $\bar{X}^0=X^0$; $\bar{X}^k/\bar{X}^0=f^k\left(X^1/X^0, ...,X^4/X^0\right)$ in $\mathfrak{H}_5$.

Having established this isomorphism of groups, between $\mathfrak{I}$ and $\mathfrak{H}_5$, Jordan "expects" that the Einstein-Maxwell theory is in a  "much more symmetric and at the same time simpler form" if described in five homogeneous coordinates $X^\mu$. This line of reasoning, called projective theory of relativity, goes back to work by Oswald Veblen and is discussed in more detail in subsection \ref{ss.veblen}. 

It begins \citep[§23]{Jordan1952} with the introduction of homogeneous coordinates:
\begin{quote}
	Let us consider the four world coordinates $x^k$ (the notation $x^{(k)}$ will now be abolished) as four independent, but otherwise arbitrary homogeneous functions of zero degree in the $X^\mu$.
	%Die vier Weltkoordinaten $x^k$ (die Bezeichnung $x^{(k)}$ soll jetzt abgeschafft werden), denken wir uns also als vier unabhängige, aber sonst beliebige homogene Funktionen nullten Grades in den $X^\mu$.
\end{quote}
That is, $x^k_{\ \ |\mu}X^\mu=0$. The task, at this point, is to define an anologous formalism of Riemannian geometry on this (five dimensional) space of which the $X^\mu$ are the coordinates of, in a way that it is induced by the known four dimensional formalism. This means that instead of considering the full group of coordinate transformations, the only legitimate invariance is the one with respect to $\mathfrak{H}_5$, i.e. with respect to coordinate transformation that are of the form $\bar{X}^\mu_{\ \ |\nu}X^\nu=\bar{X}^\mu$, as in equation \ref{eq.trafo}. The notion of a tensor is thus inappropriate, instead, objects called "projectors" are being introduced. They ensure the possibility of conversion between the five- and the four-dimensional formalism and are defined in the following way: 

\emph{The quantities $P^{\nu_1...\nu_n}_{\mu_1...\mu_m}$ are components of a projector if they transform like tensors under transformations as in equation \ref{eq.trafo} and are homogenous functions of the $X^\mu$ of degree $n-m$:}
\begin{equation}\label{eq.Euler}
	P^{\nu_1...\nu_n}_{\mu_1...\mu_m|\lambda} X^\lambda=(n-m)P^{\nu_1...\nu_n}_{\mu_1...\mu_m}.
\end{equation}
The second condition, on homogeneity, determines that knowledge of a projector field means nothing more than the knowledge of a field in a fourdimensional manifold. 

An important class of projectors are the coordinate functions $X^\mu$, which according to equation \ref{eq.trafo} transform like vectors under $\mathfrak{H}_5$. Jordan continues to call projectors $a^\nu$ vectors. Finally, in order to meaningfully translate five-dimensional projector equations into four-dimensional tensor equations, the notion of a reduction (\emph{Verkürzungen}, not to be confused with contractions) of vectors and other projectors is very important \citep[p.116]{Jordan1952}. In the case of a vector the definitions reads as $a^k=x^k_{\ |\mu}a^\mu$ - the result of a reduction of a five-dimensional vector is thus a four-vector. This procedure can be carried over analogously to arbitrary projectors.

\subsection{Covariant Derivative of Projectors and \emph{Kongruenzdifferentiation}}
%-----------------------------------------------------------------------------------
To be able to formulate differential projector equations, it is necessary to have a well-defined notion of covariant derivative on this new space. But \citep[§24]{Jordan1952}:
\begin{quote}
	The formulation of invariant differential relations between projectors does not pose a new problem for us. The concepts and theorems on covariant differentiation developed in Chapter I are also applicable here [...]. \\
	With regard to the \emph{definition} of the covariant derivative, it is best to adhere to the axiomatic characterization of this operation as described in § 9 [...].
	%Die Formulierung invarianter Differentialbeziehungen zwischen Projektoren stellt uns kein neues Problem. Die in Kapitel I entwickelten Begriffe und Sätze über kovariante Differentiation sind auch jetzt anwendbar [...]. \\
%	Betreffs der \emph{Definition} der kovarianten Ableitung halten wir uns jetzt am besten an die in § 9 ausgeführte axiomatische Kennzeichnung dieser Operation [...].
\end{quote}
 
The axiomatic definition is of advantage here as homogenous coordinates cannot be plane anywhere because of the (Euler) condition $g_{\mu\nu|\sigma}X^\sigma=-2g_{\mu\nu}$, equation \ref{eq.Euler}. Despite of the straight forward take-over of the axiomatic definition, the covariant derivative of a projector does have peculiar properties that come from the Euler condition. For example, it holds that $X_{\sigma||\rho}+X_{\rho||\sigma}=0$, that is, the Killing equation holds true for all (projective) vectors. And for the (projective) Riemann tensor $R^\nu_{\ \cdot\ \mu\sigma\tau}$, which is defined analogously to the non-projektive situation, and its contraction, the Ricci tensor $R_{\mu\tau}$, it holds $R^\nu_{\ \cdot\ \mu\sigma\tau}X_\nu=X_{\sigma||\tau||\mu}$ and $R_{\nu\tau}X^\nu=\frac{1}{2}X^\mu_{\ \cdot\ \tau||\mu}$. 

Comparing covariantly differentiated objects of the five and the four-dimensional space hits the difficulty that even if the reduction (\emph{Verjüngung}) of two projectors yields the same tensor, Jordan calls those projectors \emph{congruent},  $P^{\nu_1...\nu_n}_{\mu_1...\mu_m}\equiv Q^{\nu_1...\nu_n}_{\mu_1...\mu_m}$, the reduction of their respective covariant derivatives need not be the same tensor \citep[p.121]{Jordan1952}. For this reason, yet another "differentiationsähnliche Operation" is introduced, which Jordan calls "\emph{Kongruenzdifferentiation}", denoted by three bars, $|||$. It is defined as one of the two, then equivalent, equations
\begin{equation}\begin{aligned}
	a_{\mu|||\lambda} &= a_{\mu||\lambda}+\frac{X_{\mu\lambda}X^\sigma-X^\sigma_{\ \cdot\lambda}X_\mu}{2J}a_\sigma; \\
	a^\mu_{\ \ |||\lambda} &= a^\mu_{\ ||\lambda}+\frac{X^\mu_{\ \cdot\lambda}X_\sigma-X_{\sigma\lambda}X^\mu}{2J}a^\sigma,
	\end{aligned}\end{equation}
where $X_{\rho\sigma}:=2X\sigma||\rho=X_{\sigma|\rho}-X_{\rho|\sigma}$. An important quantity in these formulae is the "invariant", i.e. scalar field, $J$, which is defined as $J=g_{\mu\nu}X^\mu X^\nu$. It plays a decisive role in Jordan's extended gravitational theory, as is briefly explained in subsection \ref{ss.physics}.

This conjunct definition above, in any case, is necessary in order to satisfy axiom III. In order to satisfy axiom IV, it is taken as the definition for the Kongruenzdifferentiation of scalars. With this, this operation is well defined and satisfies all axioms of the covariant derivative - despite axiom V. This notion enables an important property:\\
\emph{The reduction (Verkürzung) of the Kongruenzableitung $P^{\nu_1...\nu_n}_{\mu_1...\mu_m|||\lambda}$ is the covariant derivative of $P^{j_1...j_n}_{k_1...k_m||l}$ of the reduction $P^{j_1...j_n}_{k_1...k_m}$ of $P^{\nu_1...\nu_n}_{\mu_1...\mu_m}$.}

That this theorem can easily be proven true, significantly relies on the usage of an axiomatic definition of the covariant derivative. Jordan addresses the problem of translating projector equations into fourdimensional tensor equations in the following way \citep[p.123]{Jordan1952}: 
\begin{quote}
	In its treatment, it pays off that we took the trouble in §9 to give an axiomatic characterization of the covariant derivative, which makes its definition independent of the concept of a coordinate system that is flat at a point $P$. This saves us from having to conduct a differential geometric investigation of the relationship between five-dimensional and four-dimensional differentiation, and the “algebraization” of the problem represents a radical simplification compared to the laborious proofs found in older literature [...].
	%Bei seiner Behandlung macht es sich bezahlt, daß wir uns in §9 die Mühe gemacht haben, eine axiomatische Kennzeichnung der kovarianten Ableitung zu geben, durch welche ihre Definition unabhängig wird vom Begriff des in einem Punkte $P$ ebenen Koordinatensystems. Hierdurch nämlich wird uns eine differentialgeometrische Untersuchung des Zusammenhangs fünfdimensionaler und vierdimensionaler Differentiation erspart; und die "Algebraisierung" des Problems bedeutet eine durchgreifende Vereinfachung gegenüber den mühsamen Beweisen der älteren Literatur [...].
\end{quote}

With the help of this theorem, the respective Riemannian tensors and their derived tensors can be compared. Turns out that the four-dimensional Riemann tensor equals the reduction of the projector $R^\nu_{\ \cdot\ \mu\sigma\tau} - \frac{1}{4J}\left[X_{\mu\tau} X^\nu_{\ \cdot\ \sigma}-X_{\mu\sigma} X^\nu_{\ \cdot\ \tau}+2X_{\sigma\tau} X^\nu_{\ \cdot\ \mu}\right]$. There are similar formulas for $R_{\mu\nu}$ and $R$, 
\begin{equation}\begin{aligned}
		G_{kl} &= R_{kl}-\frac{1}{2J}X_{kj}X_l^{\cdot\ j}-\frac{J_{|k||l}}{2 J}\frac{J_{|k}J_{|l}}{4 J^2}\\
		R &= G+\frac{1}{4J}X_{kl}X^{kl}+\frac{1}{J}J^{|k}_{\ \ ||k}-\frac{1}{2J^2}J^{|k}J_{|k},
	\end{aligned}
\end{equation}
where $G_{kl}$ is the four-dimensional Ricci tensor and $G=g^{kl}G_{kl}$.

\subsection{Projective Relativity and Jordan's Extension}\label{ss.physics}
%--------------------------------------------------------
The most crucial assumption within the projective theory of relativity is $J=1$. This is necessary to reduce the total number of independent field components from 15 to 14, so as to match the number of independents field components within the Einstein-Maxwell theory. Jordan emphasized \citep[p.123,128]{Jordan1952} that previous treatment of the mathematical formalism, in particular the comparison of the projector to the tensorial case, had already built in this condition. In this sense, his algebraic way to construct the formalism, which does not need this condition at this early point, is much more general. With this condition, the equations found above simplify considerably, $R_{kl}=G_{kl}+\frac{1}{2}X_{kj}X_l^{\cdot\ j}$, $R=G+\frac{1}{4}X_{kl}X^{kl}$.

The invariant $J$ is also fundamental for another reason. This is for the following property:
\begin{equation}
\left\{\left(X_{kl}/J\right)_{||j}\right\}_{[klj]}=0,	
\end{equation}
where $\left\{A_{klj}\right\}_{[klj]}$ is Jordan's notation for $A_{klj}+A_{ljk}+A_{jkl}$. This property is the basis for the following interpretation \citep[p.124]{Jordan1952}: 
\begin{quote}
	The theory thus yields, from the geometry of the five-dimensional homogeneous manifold, a \emph{four-dimensional six-vector}
	\begin{equation}
		X_{kl}/J=-X_{lk}/J,
	\end{equation}
	which has the property (12)\footnote{This number refers to our numbering of equations. In Jordan's textbook \citep{Jordan1952} this would be equation (4). Also, “§7” refers to Jordan's textbook as before.} , i.e., according to §7, it is the \emph{rotation of a four-vector}.
	%Die Theorie ergibt also von selber - aus der Geometrie der fünf-dimensionalen homogenen Mannigfaltigkeit heraus - einen \emph{vier-dimensionalen Sechservektor}
%	\begin{equation}
%		X_{kl}/J=-X_{lk}/J,
%	\end{equation}
%	der die Eigenschaft (12)\footnote{This number refers to our numbering of equations. In Jordan's textbook \citep{Jordan1952} this would be equation (4). Also, "§7" refers to Jordan's textbook as before.} besitzt, also nach §7 \emph{Rotation eines Vierervektors ist}.
\end{quote}
This now offers the interpretation of $X_{kl}$ as being proportional to the electromagnetic six-vector $F_{kl}$, $X_{kl}=\frac{\sqrt{2\chi}}{c}F_{kl}$.

The field equations are derived extremizing the integral
\begin{equation}
	\int\left(R-\lambda\left[1-J\right]\right)\sqrt{-\overset{\text{5}}{g}}dX^0...dX^4, 
\end{equation}
with the constraint $J=1$. In this formula, $\lambda$ is a Lagrange factor and $\overset{\text{5}}{g}=Det|g_{\mu\nu}|$. The corresponding Lagrange equations are
\begin{equation}
	R^{\mu\nu}-\frac{1}{2}g^{\mu\nu}R+\lambda X^\mu X^\nu=0.
\end{equation}
Going over to four dimensions, with the above formulae, one can see that this system is indeed equivalent to Einstein-Maxwell, equation \ref{eq.EM}\footnote{$\lambda$ conveniently eliminates itself.}.

Jordan concludes \citep[p.130]{Jordan1952}, 
\begin{quote}
	So we have actually transformed the Einstein-Maxwell theory into a five-dimensional symmetrical form without changing the content of its statements in any way. [...] But the five-dimensional formulation as such is much simpler than the four-dimensional one—there can be no doubt about that; [...] so we can confidently assert that it is only with the projective theory of relativity of Kaluza-Klein-Veblen that we have correctly understood the Einstein-Maxwell theory in its inner mathematical structure and gained a deeper insight into the harmonies of the laws of physics.
	%Wir haben also tatsächlich die \emph{Einstein-Maxwell}sche Theorie - \emph{ohne} den Inhalt ihrer Aussagen irgendwie zu \emph{verändern} - in eine fünfdimensional-symmetrische Gestalt gebracht. [...] Aber die fünfdimensionale Formulierung als solche ist \emph{viel einfacher als die vierdimensionale} - daran kann gar kein Zweifel bestehen; [...] so dürfen wir zuversichtlich behaupten, daß wir erst mit der projektiven Relativitätstheorie von \emph{Kaluza-Klein-Veblen} die \emph{Einstein-Maxwell}sche Theorie in ihrer \emph{innern mathematischen Struktur richtig verstanden} und einen vertieften Einblick in die Harmonien der physikalischen Gesetze gewonnen haben. 
\end{quote}
In this sense, Jordan also differentiates these efforts from other unified field theories, which, according to Jordan, mostly "want to change and generalize" the Einstein-Maxwell theory such that the latter is only a limiting case of the new theory, which is not the case here. His \emph{erweiterte Gravitationstheorie}, however, still wants to achieve exacly this:
\begin{quote}
	However, the final result of this analysis certainly compels us to consider a certain generalization of the theory. The condition [J=1] is an unattractive element of the theory and, in a sense, represents a mutilation: it seems reasonable to examine what would happen if we \emph{refrained} from using this assumption [J=const] to force the projective theory into complete agreement with the \emph{Einstein-Maxwell} theory.
	%Aber das Endergebnis dieser Analyse \emph{drängt} nun freilich doch dazu, eine bestimmte Verallgemeinerung der Theorie zu erwägen. Die Nebenbedingung [J=1] ist ein unschönes Element der Theorie, und bedeutet gewissermaßen eine Verstümmelung: Es liegt nahe, zu prüfen, was sich ergeben würde, wenn wir es \emph{unterlassen}, durch diese Annahme [J=const] gewaltsam die projektive Theorie zur vollständigen Übereinstimmung mit der \emph{Einstein-Maxwell}schen zu zwingen.
\end{quote}

For Jordan, dropping the condition $J=const$ is "equivalent to assuming that the gravitational “constant” $\chi$ is in fact a variable, a \emph{scalar field quantity}", following Dirac's idea mentioned in the introduction. And this is now of course the reason for his general, algebraic treatment of the mathematical formalism, including the axiomatic definition of the covariant derivative, with a non specified quantitiy $J$. Discussing which variational principle would be the most suitable for this new theory, Jordan settles on 
\begin{equation}
	\delta\int J^\alpha\left(R-\lambda\frac{J^{|\mu}J_{|\mu}}{J^2}\right)\sqrt{-\overset{\text{5}}{g}}dX^0...dX^4=0.
\end{equation}
Choosing $\alpha=1/2$, this is equivalent to %(and here, Jordan refers to Ludwig's paper [find paper]) 
\begin{equation}
	\delta\int\chi\left(G+\frac{\chi}{2c^2}F_{kl}F^{kl}-\left(\lambda+\frac{1}{2}\right)\frac{\chi^{|k}\chi_{|k}}{\chi^2}\right)\sqrt{-\overset{\text{4}}{g}}\ dx^1...dx^4=0.
\end{equation}
Here, $F_{kj}=CX_{kj}/J$, $C$ an arbitrary constant factor with the "right dimension" and finally, with $Jc^2/2C^2=\chi$, the variable gravitational invariant $\chi$. The variation is supposed to be with respect to all 15 field quantities $g_{kl}$, $\Phi_k$ (the four potential) and $\chi$. 

The discussion as to why this is the physically correct variational principle, Jordan postpones to the next chapter, which, in particular, does not presuppose any of the mathematical groundwork from which Jordan draws the legitimacy, in fact necessity, of his extened gravitational theory. Since the focus of this paper, in any case, is not the discussion of the physical content of Jordan's theory, a further analysis of physical considerations is omitted here.

\section{Covariant Derivative - Parallel Transport - Connection: A very brief history}\label{s.briefhistory}
%-------------------------------------------------------------------------------------
%-------------------------------------------------------------------------------------
%	
\paragraph{Covariant Derivative}
The formula for covariant differentiation first came up in Elwin Bruno Christoffel's work in the context of differential forms in 1869 \citep{Christoffel1869}, though not by this name. Like within the theory of algebraic invariants, which was developed earlier, the aim of Christoffel's work at that time was to find a calculus to treat differential invariants, i.e. to answer the question under which conditions two differential expressions can be transformed into each other and how the transformation relations would then look like \citep[p.58]{Reich1994}. Investigating these conditions for differential forms of second order, $F=\sum\omega_{ik}\partial x_i\delta x_k$, Christoffel conveniently introduced what would later be called Christoffel symbols of the first, $\left[\begin{array}{l}
	\text{$g\ h$} \\[4pt]
	\text{$\ k$}
\end{array}\right]$, and of the second kind, $\left\{\begin{array}{l}
\text{$i\; l$} \\[5pt]
\text{$\ r$}
\end{array}\right\}$, and determined their transformations properties. Based on this, he analysed how to get from transformation relations of a form of a certain order to the transformation relations of a form of one order higer. And the resulting formula is the covariant differentiation, 
\begin{equation}
	(ii_1...i_\mu)=\frac{\partial(i_1i_2...i_\mu)}{\partial x_i}-\sum_\lambda\left[ \left\{\begin{array}{l}
		\text{$i\; i_1$} \\[5pt]
		\text{$\ \lambda$}
	\end{array}\right\}
	(\lambda i_1...i_\mu)+
\left\{\begin{array}{l}
	\text{$i\; i_2$} \\[5pt]
	\text{$\ \lambda$}
\end{array}\right\}( i_1 \lambda...i_\mu)+...\right],
\end{equation}
where $(i_1i_2...i_\mu)$ denotes the components of the original form.  It is interesting to stress here that Christoffel did not think in terms of coordinates and coordinate transformations, but that it was really mainly about transforming different forms into each other.

From the late 1880s on, Christoffel's calculus was picked up by Gregorio Ricci-Curbastro who made it the center piece of what he would later call absolute differential calculus, which in turn is the basis for today's tensor calculus. And it is in Ricci's work that the term "covariant derivative" first came up \citep[p.72ff]{Reich1994}. Together with his student Tullio Levi-Civita, Ricci wrote a major treatise on his reasearch, submitted 1988 \citep{Ricci1901}.\footnote{We owe a very detailed and well written account of this part of the history of maths to historian Karin Reich in \citep{Reich1994}. Marco Giovanelli gave an interesting philosophical account in \citep{Giovanelli2012}.} Both made an effort to make their work known to and applied by physicists and some indeed picked it up\footnote{Karin Reich mentions Max Abraham and O. Tedore as examples \citep{Reich1994}.}. Albert Einstein only became aware of it in 1912 through his friend and colleague in Zürich, the geometer Marcel Grossmann, when searching for a generally covariant theory of gravitation.\footnote{Jürgen Renn, Michel Janssen and others have worked out in detail the origin story of Einstein's field equations, see for example \citep{Renn2007a}.} And it is Grossman who introduces the word "Christoffel'sche Drei-Indizes-Symbole", i.e. Christoffel symbols, while he calls the covariant derivative "Erweiterung", as it raises the order of a tensor.\footnote{Einstein in turn famously introduces the summmation convention and the notation $\Gamma_{\mu\nu}^\tau$ for the Christoffel symbol of the second kind (as components of the gravitational field).}

Even though Grossmann was a geometer, differential geometry and the absolute differential calculus was generally practised by different groups of mathematicians \citep[p.184]{Reich1994}.\footnote{There is a recent article by Alberto Cogliati \citep{Cogliati2024} which aims to "challenge this historiographical tenet".} This is, in particular, the historical origin of why Einstein never thought of general relativity as geometrizing gravity\footnote{For a more philosophical analysis, see \citep{Lehmkuhl2014}.}: he employed Ricci's absolute differential calculus in the way Ricci had thought of it - as a purely formal calculus without geometric interpretation \citep[p.210]{Reich1994}\citep[p.1]{Reich1992}\cite[p.263]{Scholz1999}. 

\paragraph{Parallel Transport}
The link between differential geometry and the Ricci calculus was finally established by solving a conceptual problem within Riemannian geometry, which was that it did not have a notion of parallelism. Studying how directions at different points could be considered parallel, Tullio Levi-Civita rediscovered in 1916 a set of differential equations previously found by other authors, but that he could now attribute a geometrical meaning to \citep{LeviCivita1917}\citep[p.81]{Reich1992}: 
\begin{equation}\label{eq.parallel}
	\frac{d\xi^{(i)}}{ds}+\sum_{n}^{1}{}_{ij}\left\{\begin{array}{l}
		\text{$j\; l$} \\[5pt]
		\text{$\ i$}
	\end{array}\right\}
	\frac{dx_j}{ds}\xi^{(l)}=0, 
\end{equation} 
where again $\left\{\begin{array}{l}
	\text{$j\; l$} \\[5pt]
	\text{$\ i$}
\end{array}\right\}=\Gamma^i_{jl}$ are the Christoffel symbols and $\xi^{(i)}, \xi^{(i)}+d\xi^{(i)}$ are the directions at the two neighbouring points which differ by (path dependent) $dx_j$.\footnote{Reich emphasizes that Gerhard Hessenberg was actually a few months earlier with very similar observations \citep{Reich1992}, wanting to introduce a more geometrically minded calculus into general relativity. However, Levi-Civita's parallel transport was somewhat more general, not only restricted to geodesics, as in Hessenberg's case. Jan Arnoldus Schouten tried to simplify the calculus underlying general relativity, thinking of tensors as quantities instead of in terms of its components. He introduced the concept of geodesically co-moving coordinate systems with respect to which the covariant differentiation is just the ordinary differentiation. With this concept, Schouten also arrived at a notion of parallelism, just somewhat later than Levi-Civita \citep[p.83ff]{Reich1992}.} This also gave a geometric meaning to the covariant derivative. In fact, it was Einstein's general relativity which motivated Levi-Civita to search for a geometrical interpretation of the formalism involved \citep[p.118]{Goodstein2018} \citep[p.263]{Scholz1999}. %Connecting to this work, Levi-Civita's colleague at Padua, Francesco Severi, went the other way round. He started from a construction of displacement along geodesics and was then able to derive equation \ref{eq.parallel} \citep[p.88]{Reich1992}.
Or, in Hermann Weyl's words \citep[p.537]{Weyl1949}, "The great importance which Riemannian geometry acquired for Einstein's theory of gravitation gave the impetus to develop this geometry further, to study more carefully its foundations and, as a consequence of such analysis, to generalize it in various directions. The first and decisive step was Levi-Civita's discovery of the notion of infinitesimal parallel vector displacement."

\paragraph{Connection}
An even deeper understanding of the covariant derivative was brought about by Hermann Weyl. His paper "Gravitation und Elektriziät" from 1918 \citep{Weyl1918} was a synthesis of mathematical, physical and philosophical ideas, which laid the basis for a second and decisive step in this conceptual lineage from covariant derivative to parallel transport - the connection. At the core of it was the observation that (infinitesimal) parallel transport made manifest a residual element of distant geometry, which he saw overthrown by general relativity: While comparing vectors at different points was no longer possible by default within Riemannian geometry, it is still possible to compare the length of vectors at two different points. Physically interpreted, this would mean the involvement of an instantanous action-at-a-discance measurement process, which Weyl saw prohibited within relativistic field theories. Instead, Weyl thought, the length of a vector should equally not be preserved during parallel transport, making the geometry "purely infinitesimal".\footnote{Historian of math Erhard Scholz has worked extensively on Hermann Weyl. See for example \citep{Scholz1994}, \citep{Scholz1995}, \citep{Scholz1999} and \citep{Scholz2001}.} Or, in Weyl's words \citep[p.538]{Weyl1949}, "Thus an affine infinitesimal geometry has sprung up beside Riemann's metric one."

In a subsequent article on "Reine Infinitesimalgeometrie" \citep{Weyl1918b}, Weyl defined the notion of an "affine connection" (\emph{"Affiner Zusammenhang"}) in the following way \citep[p.389]{Weyl1918b} \citep[p.90]{Reich1992}:
"Ist P' ein zu dem festen Punkt P unendlich benachbarter, so \emph{hängt} P' mit P \emph{affin zusammen}, wenn on jedem Vekor in P feststeht, in welchen Vektor in P' er durch \emph{Parallelverschiebung} yon P nach P' übergeht." The adjective \emph{affine} comes from the condition that  "Die Verpflanzung der Gesamtheit der Vektoren von P nach dem unendlich benachbarten Punkte P' durch Parallelverschiebung liefert eine affine Abbildung der Vektoren in P auf die Vektoren in P'." Weyl further defines the $\Gamma^i_{rs}=\Gamma^i_{sr}$ to be the "components", i.e. the central quantities, of the affine connection. With the affine connection he laid out in \citep{Weyl1918}, he famously provided the first unified field theory.\footnote{Which at the same time was also the first so-called gauge theory; see \citep{ORaifeartaigh1997} for more on that angle}

There were several approaches that tried to generalize Weyl's notion of an affine connection further \citep[p.96]{Reich1992}. The most general and thorough treatment of this was probably done by Jan Arnoldus Schouten, who preferred the word \emph{Übertragung} over Weyl's \emph{Zusammenhang}. In his 1922 paper "Über die verschiedenen Arten der Übertragung in einer $n$-dimensionalen Mannigfaltigkeit, die einer Differentialgeometrie zugrunde gelegt werden können" \citep{Schouten1922} he explained, 
\begin{quote}
	The work of the authors mentioned above has gradually shown that this Übertragung [...] can also be defined in a much freer way and even completely independently of a fundamental tensor. However, there are still a few steps to be taken before we can arrive at a complete overview of the various possibilities offered by extended differential geometry.
	%In den Arbeiten der erwähnten Autoren hat sich nun nach und nach gezeigt, daß diese Ubertragung [...] auch in viel freierer Weise und sogar ganz unabhängig von	einem Fundamentaltensor definiert werden kann. Es fehlen aber noch einige Schritte, um zu einer vollständigen Übersicht der verschiedenen Mögichkeiten der erweiterten Differentialgeometrie zu gelangen.
\end{quote}
In this paper, Schouten then sets out to find (and did find) the most general \emph{lineare Übertragung}\footnote{I keep the word Übertragung also in the English translation of the quotes.}, which impressively gives way to 18 different cases \citep[p.73]{Schouten1922}. 

It is instructive to take a closer look at how Schouten sets out to find this most general Übertragung. He begins with the following observation \citep[p.63]{Schouten1922},
\begin{quote}
	In principle, one would even be completely free in the choice of Übertragung. For example, one could specify any Übertragung for each size at each point of $X_n$ and for each direction. However, the corresponding differentiation would then not satisfy any of the formal laws to which ordinary differentiation is subject.
	%Prinzipiell wäre man in der Wahl der Übertragung sogar vollständig frei. Man könnte z.B. für jede Größe in jedem Punkt der $X_n$, und für jede Richtung eine ganz beliebige Übertragung festlegen. Die zugehörige Differentiation würde dann aber keinem der formalen Gesetze genügen, denen die gewöhnliehe Differentiation unterworfen ist.
\end{quote}
Schouten therefore gives a list of properties that the covariant derivative $\delta$ resulting from an (yet to be constructed) Übertragung needs to obey:
\begin{enumerate}
	\item A quantity and its differential have the same number of components, and these transform in the same way.
	%Eine Größe und ihr Differential haben dieselbe Anzahl Bestimmungszahlen (components), und diese transformieren sich in derselben Weise.
	\item If $\Phi$ is any quantity of any degree, then $$\delta\Phi=\frac{\delta\Phi}{\delta x^\mu}dx^\mu.$$
	\item If $\Phi$ and $\Psi$ are arbitrary quantities of arbitrary degree, then $$\delta(\Phi+\Psi)=\delta\Phi+\delta\Psi.$$
	\item Similarly, $$\delta(\Phi\Psi)=\delta\Phi\Psi+\Phi\delta\Psi.$$
	\item The differential of a scalar is identical to the ordinary differential.
\end{enumerate}

An \emph{Übertragung} whose covariant derivative fullfills these properties is then called \emph{lineare Übertragung}. While a characterization of the covariant derivative is given here, in terms of preferred properties, it is important to note that this is not an axiomatic definition of the covariant derivative (like Jordan has provided). Schouten is searching here for a connection whose covariant differentiation has these (characteristic) properties. He does not provide an analysis as to whether these properties comprise a complete axiomatic basis for the covariant derivative, nor is it here a primary and independent notion.

\paragraph{The Covariant Derivative within General Relativity and Unified Field Theory}
For work on general relativity proper, that is, for work within the paradigm that Einstein set out in 1915/1916 with \citep{Einstein1915b} and \citep{Einstein1916a}, there was no reason and room for altering the mathematical formalism, including the covariant derivative.

For someone who was interested in developing general relativity further, along the lines of unified field theory, testing out other differential geometrical structures was the preferred, if not the only, route. Put the other way round, it is probably safe to say that the avenue that was opened up by the concept of connection enabled the unified field theory program in the first place.\footnote{There were long standing physical motivations as well, of course. However, if this conceptual degree of freedom would not have opened up, it seems entirely unclear how these would have manifested formally (and if at all).} 
Instead of listing the manifold ways in which unified field theory has devised new affine geometries, we refer here to already existing historical literatur which has got it covered already: \citep{Goenner2004}, \citep{Goenner2014}, \citep{Vizgin2011}, \citep{Dongen2010}, \citep{Sauer2014} and others.

\section{Projective Geometry and Projective Relativity}\label{s.projective}
%-----------------------------------------------------------
%-----------------------------------------------------------
Following the developments that lead to the notions of parallel transport and connection, attempts were made to incorporate projective geometry into an affine framework.\footnote{The classic construction of projective spaces is to "add" for each direction a point in infinity ("projective closure"), at which lines meet and are thus no longer parallel.} Above all by Weyl himself, who pioneeringly wrote in 1921 \citep[p.99]{Weyl1921c}, "\emph{Projective} and \emph{conformal} geometry arise through abstraction from affine and metric geometry, respectively." Also Schouten in his book "Der Ricci Kalkül" vom 1924 \citep{Schouten1924} treats "projective curvature" (\emph{Projektivkrümmung}) as part of (his) affine geometry.
However, as also Weyl later acknowledged \citep[p.538]{Weyl1949}, referring to work by Luther Pfahl Eisenhart and Oswald Veblen (Princeton), Schouten and David van Dantzig (Delft), Elie Cartan and Shiing Shen Chern (Paris) and others: 
\begin{quote}
In several ways these authors soon arrived at the conclusion that it is better to establish projective differential geometry not by abstraction from the affine brand, as described above, but independently, namely by associating with each point $P$ of the manifold a projective space $\Sigma_p$ in the sense of Poncelet and Plücker, this homogeneous space taking the place of the affine vector compass in the affinely connected manifold.
\end{quote}

These three schools mentioned by Weyl, however, started out with rather different viewpoints on this matter. While Schouten and Cartan tried to stay true to Felix Klein's Erlangen Programm, regarding geometry as the invariance theory of a (Lie) group (see. e.g. \citep{Schouten1926} and \citep{Cartan1928})\footnote{For more on Klein's programme see \citep{Rowe2025}.}, Veblen and Eisenhart deviated completely from this and put forward a notion of geometry based on systems of paths. Interestingly, a projective account of unified field theory stood at the end of both of these strands.\footnote{Of course, Schouten's and Cartan's approaches can be distinguished further and it was only Schouten and his school who came up with a projective theory of relativity, however we will not go into this level of detail here.}  In this section, however, only the development of the geometric program of Veblen-Eisenhart which, I believe, portrays best the struggles geometers went trough in the 1920s. It is also a less straightforward approach to a new concept of geometry, contrasting nicely with the programmes of Schouten and Cartan.
%Cartan saw his generalization of geometry as being in the spirit of Felix Klein's \emph{Erlanger Program}, only based on Lie Groups, rather that the classical groups that Klein was thinking of as the basis for geometry. 

%Form \citep[p.182]{Veblen1929a}
%\begin{quote}%
	%Nevertheless the hold of the Erlanger Programm upon the imagination of
	%mathematicians is such that attempts were sure to be made to revamp the Programm
	%so as to adapt it to the new order of things. And these attempts have
	%had a considerable degree of success. The concept of infinitesimal parallelism
	%which had been introduced by LEVI-CIVITA was developed and enlarged by WEYL
	%and has been generalized by a number of mathematicians. In particular, CARTAN
	%and SCHOUTEN have shown that there are other ways than those forseen by
	%KLEIN of connecting up the theory of continuous groups with geometry. As
	%CARTAN has said, we may regard a Riemannian space as a non-holonomic
	%Euclidean space, and many of the generalizations of Riemannian spaces can be
	%arrived at in a similar manner.
%\end{quote}

%--------------------------
\subsection{The Geometry of Paths Program}\label{ss.veblen}
%-----------------------------------------------------------
Oswald Veblen had started his career, before the First World War, with work on axiomatic foundations of (projective) geometry (for example, \citep{Veblen1904} and \citep{Veblen1908}). It was only later that he underwent a "differential-geometric turn" \citep[p.11]{Ritter2011}, which involved a strong interest in general relativity and unified field theory. His colleague at Princeton, Luther Pfahl Eisenhart, was already America's expert on the theory and the one who invited Einstein to come lecture in Princeton in 1921.\footnote{These lectures would turn into Einstein's only textbook on his theory, the famous "The Meaning of Relativity" \citep{Einstein1922}.}

In the early 1920s, Veblen and Eisenhart laid out an ambitions, school-building, reseach program that integrated physics and mathematics and which would last for about a decade - the \emph{geometry of paths program}. Even though intellectually rooted in Weyl's theory of connections and unified field theory, they made a point in not wanting to generalize geometry (and advance physics) further in this direction. Instead, they put primary the motion of particles, that is, the paths particles take in physical circumstances. In the case of general relativity, they take the shape of geodesics,\footnote{In original notation of e.g. \citep{Veblen1923a} and \citep{Eisenhart1922}.} 
\begin{equation}\label{eq.geodesic}
	\frac{d^2x^i}{ds^2}+\sum_{\alpha,\beta}\Gamma^i_{\alpha\beta}\frac{dx^\alpha}{da}\frac{dx^\beta}{ds}=0. 
\end{equation}
In an affine setting, the $\Gamma$ are determined by the connection. Veblen and Eisenhart, however, now propose to drop this condition and look at systems of paths that are determined by a set of differential equations like in \ref{eq.geodesic} where the $\Gamma$ are arbitrary functions (albeit symmetric in the lower indices).

This approach, they saw as "a generalization of both of the earliest part of elementary geometry and of some of the most refined of physical theories" \cite[p.137]{Veblen1923a}, thus connecting Veblen's earlier work on the axiomatic foundation of geometry with the more recent developments in differential geometry, including general relativity.\footnote{A detailed philosophical account of how these things go together is given in \citep{Veblen1923a}.}
Veblen further explains that \citep[p.136]{Veblen1923a}\footnote{This approach is also interesting for another reason: The way the motion of particles is derived from a field theory played a major role for Einstein (this is called "the problem of motion", see e.g. \citep{Lehmkuhl_2019}). Veblen and Eisenhard seem to avoid this problem by starting with paths.}
\begin{quote}
	The intuitive idea suggested by this name [the geometry of paths] is that we are dealing not with the empty space of analysis situs, but with a manifold in wich we find our way around by means of paths. It may also serve to remind us that we have a generalization of an inertial field, for the characteristic of a field of inertia is that through every point and in every direction, there is a path with may be taken by a free particle. 
\end{quote}

Allowing the $\Gamma$s to be arbitrary, Veblen and Eisenhart show \citep[p.20]{Eisenhart1922} that in general there exists no quadratic form for which the paths are geodesics. In particular, the invariance group of a space of paths is in general the identity. This is a point that Veblen and Eisenhart stress several times (for example also here \citep[p.137]{Veblen1923a}). It is important, because in this way, their program deviates from Felix Klein's Erlangen Program which defines geometry as the invariance theory of a group.
Still, within their theory, they can derive meaningfully a curvature tensor and "[t]he theory of covariant differentiation [...] can be generalized at once to the geometry of paths by replacing the Christoffel symbols  $\left\{\begin{array}{l}
	\text{$j\; k$} \\[0.5pt]
	\text{$\ i$}
\end{array}\right\}$ by the functions $\Gamma^i_{jk}$ in all formulas."\footnote{An elaborate treatment of this can be found in \citep{Eisenhart1927}. This seems to be the last major output of the group before emphasizing the usage of invariants over paths.}

The geometry of paths program also had an explicit projective-geometric part. At first, the \emph{projective geometry of paths} was regarded as a special case, analogous to Weyl's treatment \citep[p.136-37]{Veblen1923}: 
\begin{quote}
	From the work of Weyl\footnote{They refere here to Weyl's paper \citep{Weyl1921c} mentioned before.}, who has made the most important contributions to the geometry of paths, it follows that there is not only an affine but also a projective geometry of paths. The affine geometry consinsts of those theorems which deal with the concept of infinitesimal parallelism defined by means of the functions $\Gamma$. The projective geometry deals with those properties of the paths that are so general as to be independent of any particular definition of parallelism. That there is such a projective geometry follows from the fact that more than one set of functions $\Gamma$ can be found to define the same set of paths.\footnote{A more technical statement can be found in \citep{Veblen1923}, Section 5. Projective geometry of paths.}
\end{quote}

\subsection{Differential Invariants and Geometry}
%-------------------------------------------------
The following years saw a significant shift of Veblen's research program, which can be understood as the outline of a new general viewpoint to understanding geometry. It seems that the ingredients which enabled this shift were threefold: 
1) The awareness that Klein's Erlangen program does not fit with the geometry of paths program, 
2) The efforts to develop geometry along the lines of Klein's program by other contempory geometers such as Schouten and Élie Cartan\footnote{As, for example, expressed in a lecture by Cartan at the 1924 International Congress of Mathematicians in Toronto on "La théorie des groupes et les recherches récentes de géométrie différentielle" \citep{Cartan1928}, revoking Klein's Erlangen program with Lie groups, instead of the original classical groups (see also \citep[p.168]{Ritter2011}).} 
and 3) A commissioned monograph for the series "Cambridge Tracts in Mathematics and Mathematical Physics" on "Invariants of Quadratic Differential Forms" that appeared in 1927 \citep{Veblen1927}. 

Veblen's engagement with the theory of invariants for his monograph appeared to provide him with means to address Klein's Erlangen programme.\footnote{Though the book itself is purely formal algebraic and explicitely treats no (differential) geometry.} He presented his new vision for geometry in his lecture at the International Congress of Mathematicians in Bologna in 1928 \citep{Veblen1929a}. 
About Klein's Erlangen program Veblen said, "This point of view was the dominant one for the first half century after it was enunciated. It effectively took account of subjects like Projective Geometry
which the Riemannian point of view seemed to overlook" \citep[p.181]{Veblen1929a}. However, 
\begin{quote}
	With the advent of Relativity we became conscious that space need not be looked at only as a « locus in which », but that it may have a structure, a field-theory, of its own. This brought to attention precisely those Riemannian geometries about which the Erlanger Programm said nothing, namely those whose group is the identity. In such spaces there is essentially only one figure, namely the space structure as a whole. It became clear that in some respect	the point of view of RIEMANN was more fundamental than that of KLEIN. [...]
	
	Once we have recognized that there are geometries which are not invariant theories of groups in the simple sense which we had in mind at first, we are on the way to recognize that a space may be characterized in many other ways than by means of a group. For example, there is the fundamental class of spaces of paths studied by EISENHART and some of my other colleagues, which are characterized by the presence of a system of curves such that each pair of points is joined by one and only one curve of the system. Whether or not these spaces can be characterized in other ways there can be no doubt of the significance of this way of viewing them.
\end{quote}

To Veblen, Riemman's point of view was the situation of a topological manifold, with analytic coordinate transformations, and on it a general theory of differentials. And he goes on reffering to the developments in this line of reasoning, mentioning Lipschitz, Christoffel and Ricci, and finally the mathematical physicists who were continually developing these "extremely formal and narrow in outlook" ideas in differential invariant theory. This development, so Veblen, has now "led to a conception of a differential invariant which is well suited to the comparative study of geometries." And what he wanted to understand as a geometry from now on was "the theory of one ore more such invariants".

\subsection{Generalized Projective Geometry and Projective Relativity}
%---------------------------------------------------------------------
Veblen's differential invariant approach to projective geometry is the following. Starting from what we know, "[p]rojective geometry is the theory of the straight lines free from some of the restrictions imposed by the affine treatment. One of these restrictions is that the differential equations (4) imply a particular assignment of the parameter $t$ to the points of the line" \citep[p.185]{Veblen1929a}. Here, differential equations (4) is the geodesic equation in Euclidean space, $\frac{d^2x^i}{dt^2}=0$. If this condition on the parameter $t$ is being dropped, this equations turns into
\begin{equation}
	\left.\frac{d^2x^i}{dt^2}\middle/\frac{dx^i}{dt}\right.=\varphi\left(x,\frac{dx}{dt}\right), 
\end{equation}
with $\varphi$ an arbitrary function that is homogeneous of degree one in $\frac{dx^i}{dt}$. 

If the set of equations of the system of paths is changed in this manner, the components of the affine connections $\gamma$ change in the following way: 
\begin{equation}
	\Gamma^i_{jk}+\delta^i_j\varphi_k+\delta^i_k\varphi_j,
\end{equation}
where the $\varphi$ here do not depend on $\frac{dx^i}{dt}$. These changes, of the equation of motion and the components of the connection, do not affect, however, the following quantities which can thus be considered to be determined by the Euclidean geodesic equations: 
\begin{equation}
	\Pi^i_{jk}=\Gamma^i_{jk}-\frac{1}{n+1}\left(\Gamma^a_{aj}\delta^i_k+\Gamma^a_{ak}\delta^i_j\right).
\end{equation}
These are the components of the invariant that Veblen calls the "projective connection". A projective connection for which there exists a coordinate system in which these components vanish, describes the classical projective geometry. In this sense, the theory of this invariant may be considered as a generalized projective geometry. 

To study this generalized projective geometry further, Veblen had introduced projective tensors and related operations in a short note already before his lecture in Bologna \citep{Veblen1928}. 
A projective tensor is an invariant that transforms against coordinate change like multiplication with the quantities $u^i_j=\frac{\partial u^i}{\partial \bar{x}^j}$, $u^0_i=\frac{\partial\log u}{\partial\bar{x}^i}$, $u=|\frac{\partial x}{\partial\bar{x}}|$. A projective vector $A_\alpha$ of weight $N$, for example, transform like $\bar{A}_\alpha = u^N A_\alpha u\sigma_\alpha$.
A projective tensor has $(n+l)^k$ components instead of $n^k$ as tensors on spaces of the same dimension. 

There is also "a process of projective differentiation analogous to covariant differentiation"  \citep[p.188]{Veblen1929a}:
\begin{quote}
	In this process we use an invariant called the extended projective connection with $(n+1)^3$ components which is in a simple relationship with the original projective connection. By a suitable elimination between the law of transformations of this invariant and that of the derivatives of the components of a projective tensor we find a formula which leads from any given projective tensor to another projective tensor with one more covariant index. This is the process of projective differentation.
\end{quote}
For a covariant projective tensor of the first degree $A_\alpha$ (a projective vector) 	with weight $N$ the projective derivative is given by \citep[p.162]{Veblen1928}, 
\begin{equation}
	A_{\alpha,i}=\frac{\partial A_\alpha}{\partial x^i}-A_\lambda\Pi^\lambda_{\alpha i} \text{ and } A_{\alpha,0}= -N(n+1)A_\alpha-A_\lambda \Pi^\lambda_{\alpha0}.
\end{equation}
Veblen interprets his analysis on the generalized projective geometry în the way that "we arrive at formulas which include and clarify those obtained by the geometers who have been studying the question from the point of view of infinitesimal displacements." \citep[p.189]{Veblen1929a}

During a guest lectureship at Oxford in the academic year 1928/29, Veblen got to know about Oskar Klein's work on the five-dimensional approach to unified field theory \citep{Klein1926} based on Theodor Kaluza's pionieering paper from 1921 \citep{Kaluza1921}. The historian Jim Ritter attributes this to the influence of a young Oxford physicist, Banesh Hoffmann \citep[p.170]{Ritter2011}. When Hoffmann came to Princeton afterwards, him and Veblen worked out a reformulation of the so-called Kaluza-Klein theory that makes use of Veblen's generalized projective geometry. A publication witnessing this work appeared in 1930 under the title "Projective Relativity" \citep{Veblen1930}. The abstract reads: 
\begin{quote}
	In this paper we show that the formalism of O. Klein's version of the five-dimensional
	relativity can be interpreted as a four-dimensional theory based on projective
	instead of affine geometry. The most natural field equations for the empty spacetime
	case are a combination into a single invariant set of the gravitational and electromagnetic
	field equations of the classical relativity without modification. This seems
	to be the simplest possible solution of the unification problem.
	When we drop a restriction on the fundamental projective tensor which was imposed
	in order to reduce our theory to that of Klein a new set of field equations is obtained
	which includes a wave equation of the type already studied by various authors.
	The use of projective tensors and projective geometry in relativity theory therefore
	seems to make it possible to bring wave mechanics into the relativity scheme.
\end{quote}

Mathematically, it rests on a somewhat earlier paper written by Veblen alone \citep{Veblen1930a}\footnote{In 1933, Veblen published a lon treatise on this topic, which is written in German and based on lectures he gave previously in Göttingen, Vienna and Hamburg \citep{Veblen1933}.}, while this (short) paper fleshes out only the physics. In there, the field equations they arrive at read
\begin{equation}
	\pi_{\alpha\beta}-\Phi_\alpha\Phi_\beta\pi=0, 
\end{equation}
where $\Pi_\alpha\beta$ is the projective version of the Einstein tensor and $\Phi_\alpha$ is a projective vector that relates to the electromagnetic potential.

\section{Conclusion, or: Reassessing Jordan's "Erweiterte Gravitationstheorie" in the light of its history}\label{s.reassesing}
%--------------------------------------------------------------------------
At about the same time, also Schouten and David van Dantzig came forward with a projective reformulation of five-dimensional unified field theory \cite{Schouten1931}\citep{Schouten1932}. Besides this, Albert Einstein and Walther Mayer published a unified field theory that replaces the ordinary tangent space with five-dimensional vector spaces at each point \citep{Einstein1931a}. The physicist Wolfgang Pauli summarized the situation as follows \citep{Pauli1933}:
\begin{quote}
	Van Dantzig deserves credit for having thoroughly investigated projectors in homogeneous coordinates, their covariant differentiation using three-index symbols $\Gamma^\nu_{\lambda\mu}$, geodesic lines, and the metric introduced by the invariant form $$g_{\mu\nu} X^\mu x^\nu$$ with the fundamental tensor $g_{\mu\nu}=g_{\nu\mu}$. Finally, with the help of this general projective differential geometry, Schouten and van Dantzig have provided a formulation of field theory (we are referring here only to its classical part, which corresponds to the absence of material particles) that combines all the advantages of the formulations of Klein-Kaluza and Einstein-Mayer and avoids all their disadvantages.
	%Es ist das Verdienst von van Dantzig, die Projektoren bei homogenen Iioordinaten, ihre kovariante Differentiation mittels Dreiindizessymbolen $\Gamma^\nu_{\lambda\mu}$, die geodatischen Linien und die durch die invariante Form $$g_{\mu\nu} X^\mu x^\nu$$ eingeführte Metrik mit dem Fundamentaltensor $g_{\mu\nu}=g_{\nu\mu}$ eingehend untersucht zu haben. Schließlich haben Schouten und van Dantzig mit Hilfe dieser allgemeinen projektiven Differentialgeometrie eine Formulierung der Feldtheorie gegeben (wir sprechen hier zunachst nur von deren klassischem Teil, der dem Fall der Abwesenheit yon materiellen Teilchen entspricht), die alle Vorteile der Formulierungen von Klein-Kaluza und Einstein-Mayer verbindet, und alle deren Nachteile vermeidet.
\end{quote}

And it seems that Jordan's formulation is indeed closest to Schouten and van Dantzig (which is why we don't go into further into the details of their approach here). Pauli himself explores this topic as well, also together with his co-worker Jacques Solomon, in several publications in the early 1930s \citep{Pauli1932a}\citep{Pauli1933}\citep{Pauli1933a}. But eventually, he would drop the topic. 
In any case, none of the aforementioned researchers employed the covariant derivative in an axiomatic way, like Jordan did. In a sense, Jordan's entire approach to the topic is rather idiosyncratic. 

It is still unclear, how Jordan became aware of projective geometry in the first place. He may have picked it up from Pauli, to whom he was close. He does conscientously cite Veblen, Schouten et al in his work, but when comparing his axiomatic definition of the covariant derivative to other work, he does not think of the mathematicians at all\citep[p.34]{Jordan1952}:
\begin{quote}
	\emph{Generalizations} of \emph {Riemannian} geometry, initiated by \emph{Weyl} and carried out by \emph{Eddington} and \emph{Schrödinger}, completely abandoned the requirement of a metric $ds^2=g_{kl}dx^kdx^l$ and \emph{independently} of this \emph{postulated} an operation of the covariant (or “affine”) derivative. Axiom III is then omitted, and \emph{Schrödinger} also abandoned Axiom V (which \emph{Eddington} still retained). Intuitively, this changed view would mean that the \emph{parallel transport} (described by more general coefficients $\Gamma^m_{kj}$) is considered more fundamental than the metric (which was only defined retrospectively in those theories and treated as less important): The \emph{gyrocompass} is regarded here as an elementary instrument compared to the \emph{scale}. 
	
	However, we do not wish to discuss these theories here, although our tendency (pursued for other reasons) to completely separate the concept of covariant derivative from the considerations of §8 reveals a certain similarity to the \emph{Eddington-Schröderschen} line of thought.
	%\emph{Verallgemeinerungen} der \emph{Riemannschen} Geometrie, die durch \emph{Weyl} angebahnt und durch \emph{Eddington} und \emph{Schrödinger} ausgeführt wurden, haben die Voraussetzung einer Metrik $ds^2=g_{kl}dx^kdx^l$ ganz fallen gelassen und \emph{unabhängig} von dieser eine Operation der kovarianten (oder "affinen") Ableitung \emph{postuliert}. Es entfällt dann das Axiom III., und \emph{Schrödinger} hat auch das (von \emph{Eddington} noch beibehaltene) Axiom V. aufgegeben. Anschaulich würde diese veränderte Auffassung bedeuten, daß man die (durch verallgemeinere Koeffizienten $\Gamma^m_{kj}$ beschriebene) "\emph{Parallelverschiebung}" für grundlegender hält, als die (in jenen Theorien erst nachträglich definierte und als weniger gewichtig behandelte) Metrik: Der \emph{Kreiselkompass} ist hier als ein im Vergleich zum \emph{Maßstab} elementares Instrument angesehen. \\
	%Wir wollen jedoch auf diese Theorien hier nicht eingehen, obwohl unsere (aus anderen Gründen verfolgte) Tendenz, den Begriff der kovarianten Ableitung ganz von den Erwägungen von §8 zu lösen, eine gewisse Verwandtschaft mit den \emph{Eddington-Schröderschen} Gedankengängen ergibt.
\end{quote}

The idea to link projective relativity to his cosmological concern, however, seems to be genuinely his. As mentioned in a footnote in the introduction, Jordan's extended theory of gravitation is a mathematical treatment of an hypothesis originally due to Paul Dirac. In a 1937 Nature "Letter to the Editor", Dirac speculated that "the gravitational 'constant' must decrease with time, proportionally to $t^{-1}$" \citep{Dirac1937}. Jordan quickly picked it up and published his speculations on this in \citep{Jordan1939}, which appeared 1939. Afterwards, Jordan seemed to have dropped the topic again, however. It was only in 1944, when Jordan would get back to it. In \citep{Jordan1944} he announced that new astronomical data on star formation, that were communicated to him by Albrecht Unsöld, would make it worthwile to pick up Dirac's hypothesis again. A footnote in this paper was also the first time Jordan announces work on his extended gravitational theory (though not yet under this name). A published paper on projective relativity by Jordan would, however, appear only in 1945 \citep{Jordan1945}. An intense period of work followed, also with several co-workers, which had a preliminary peak with Jordan's textbook \citep{Jordan1952}. Which role the projective formalism took in Jordan's re-evaluation of Dirac's hypothesis and at what point he became aware of its usefullness (and how) still seems unclear.

That Jordan chose an axiomatic definition of the covariant derivative, in any case, might be a late manifestiation of his academic upbringing in Göttingen, the epicenter of axiomatics.\footnote{He did already employ axiomatic formulations for Quantum Field Theory, together with Eugene Wigner.}

\section*{Acknowledgments}
This work was funded by the European Research Council (ERC) under the European Union's Horizon 2020 research and innovation programme (grant agreement $\textnormal{n}^{\circ} 101088528$, COGY). I thank the Lichtenberg Group for History and Philosophy of Physics of the University of Bonn for discussions and inspiration.

%----------------------------------------------------------------------------------------
%	BIBLIOGRAPHY
%----------------------------------------------------------------------------------------
	
	\renewcommand{\refname}{\spacedlowsmallcaps{References}} % For modifying the bibliography heading
	
	\bibliographystyle{unsrt}
	
	\bibliography{C:/Users/Admin/Documents/Archives/Bibliothek/Bibliothek} % The file containing the bibliography
	
	%----------------------------------------------------------------------------------------
	
\end{document}